%% file: alessio_fumagalli.tex
\documentclass[graybox]{svmult}

\usepackage{mathptmx}       % selects Times Roman as basic font
\usepackage{helvet}         % selects Helvetica as sans-serif font
\usepackage{courier}        % selects Courier as typewriter font
\usepackage{type1cm}        % activate if the above 3 fonts are
\usepackage{makeidx}         % allows index generation
\usepackage{graphicx}        % standard LaTeX graphics tool
\usepackage{multicol}        % used for the two-column index
\usepackage[bottom]{footmisc}% places footnotes at page bottom
\usepackage{amsmath, amsfonts, amssymb}
\usepackage{color}
\usepackage{tikz}
\usepackage{mathtools}
\usepackage{bm}
\usepackage{ulem}

\DeclareMathOperator*{\tr}{\operatorname{tr}}
\newcommand{\abs}[1]{\vert #1 \vert}

\usepackage{lineno}
\modulolinenumbers[5]
\linenumbers

\makeindex             % used for the subject index
\begin{document}

\title*{Reactive flow in fractured porous media}
\titlerunning{Reactive flow in fractured porous media}
\author{Alessio Fumagalli and Anna Scotti}

\institute{
% Put the first author's name here
Alessio Fumagalli \at  Politecnico di Milano,\\
    via Bonardi 9, 20133 Milan, Italy \\
    \email{Alessio.Fumagalli@polimi.it}
\and % Separate the authors by \and
Anna Scotti \at Politecnico di Milano,\\
    via Bonardi 9, 20133 Milan, Italy \\
    \email{Anna.Scotti@polimi.it}
}% end institute

\maketitle

\abstract{
    In this work we present a model reduction procedure to derive a
    hybrid-dimensional framework for the mathematical modeling of reactive
    transport in fractured porous media. Fractures are essential pathways in the
    underground which allow fast circulation of the fluids present in the rock
    matrix, {often characterized by low permeability}.  However, due to
    infilling processes fractures may {change their hydraulic properties and
    become} barriers for the flow creating impervious blocks.
    The geometrical as well as the physical properties of the fractures require
    a special treatment to allow the subsequent numerical discretization to be
    affordable and accurate.  The aim of this work is to introduce {a simple yet complete
    mathematical model to account for such diagenetic effects where chemical
    reactions will occlude or empty portions of the porous media and, in
    particular, fractures.}
\keywords{Reduced modeling, hybrid-dimensional framework, fracture porous media,
reactive flow
\\[5pt]
{\bf MSC }(2010){\bf:}
% msc 2010 classification of your contribution
% the list of the msc is available here: http://www.ams.org/mathscinet/msc/
80A32, 76S05, 35Q86
}
}

%\classification{}

% ====================================
\section{Introduction}\label{sec:introduction}

Fractures  play a crucial role in determining fluid flow in a geological system.
However, two critical parameters make the modelisation of fractures challenging
from both mathematical and numerical points of view. These are their apertures,
which normally are several order of magnitude smaller than any other dimensions
in the problem, and their microscopic {structure: fractures can be open or filled by
porous materials}. Fractures thus can behave as highly conductive flow pathways
that link distant parts of the geological system and allow for fast circulation
of fluid or, on the opposite side, can be clogged preventing the flow and
creating impervious parts which are not reachable.
A fracture can have a portion of its core partially or
fully filled and another portion empty. Such complexity normally creates
problems
for classical models.

In addition the fluids present in the underground can carry ions of different
types that, under certain thermal conditions, might interact and react forming
salts that precipitate and attach to the walls of the void spaces of the porous
media.  This process tends to reduce the void spaces with a direct impact on the
flow properties of the system. {We will call } these salts ``precipitate''
while the ions are called ``solutes''. Conversely, if a precipitate is already
present and the environmental conditions are such that it can dissolve we will
have an increment of the porosity (i.e. void space) and the creation of ions
that can be transported by the liquids. Some reference on this subject are:
reactive transport on porous media at pore-scale \cite{Hornung1994,Duijn2004,Noorden2009a}
and at macro-scale \cite{Knabner1995,Emmanuel2005,Agosti2015} with
experiment comparison \cite{Katz2011}. For an micro to macro
upscaling procedure see \cite{Noorden2009,Ray2019}.

In presence of fractures the situation is even more complex. The deposition or dissolution reactions can also take place
inside the fractures, substantially altering their physical properties and
impacting the global flow properties of the geological system. This work aims to
introduce a mathematical model to accurately describe this phenomena with the
technique of dimensionality reduction. This technique is rather standard in the
treatment of problem with thin interfaces and frequently used in problems
involving fractures.  Single-phase flow
\cite{Alboin2002,Martin2005,DAngelo2011,Formaggia2012,Sandve2012,Fumagalli2017a,Antonietti,Scotti2017,Boon2018,Chave2018,Stefansson2018,Berre2019b},
two-phase flow \cite{Jaffre2011,Fumagalli2012d,Elyes2015}, passive transport
\cite{Fumagalli2012a,Chave2019,Ambartsumyan2019}, and poro-elasticity
\cite{Ganis2014,Bukac2017,Ucar2018,Berge2019} {are some of the physical problems
which have been successfully modeled with this technique}.

This work is organised as follow: in Section \ref{sec:rectiveflow} the
mathematical model for flow and reactive transport in a porous media is
presented. Section \ref{sec:reducedmodel} is devoted to the derivation of the
reduced model to describe the fracture flow and transport via reduced models. Finally,
Section \ref{sec:conclusion} contains the conclusion of this work.

\section{Reactive flow}\label{sec:rectiveflow}

In this section we present the mathematical model which describes the
flow in porous media and the transport of several species of ions (solutes). These may react forming a salt and the salt may in turn dissolve forming the ions. We consider two possible
reactions: \textit{i)} the precipitation or crystallisation where these solutes form a
solid part or \textit{ii)} a
dissolution. For simplicity in the exposition, we consider
only one precipitate.

The porous media is described by the domain $\Omega \subset \mathbb{R}^n$, for
$n=2,3$. We suppose that the porous media is saturated by a single liquid phase,
e.g. water, and the ions
are transported by its motion. Finally, the fine scale composition of the porous
media is such that a Darcy model at macroscopic scale can be applied. Our
presentation is indeed given at the macro-scale. For simplicity, the initial time is assumed
to be $0$.

\subsection{Reactive model}

Let us consider several solutes $\{U_i\}_{i=1}^{N}$ which are transported in the
porous media by a liquid phase. As already mentioned, these solutes may react to form a solid part
$W$. The integer $N$ indicates the number of species of ions that are involved in the
chemical reactions,which can be written as
\begin{align}\label{eq:chemical}
    \sum_{i=1}^{N} \alpha_i^+ U_i \leftrightarrow W + \sum_{i=1}^{N} \alpha_i^-
    U_i.
\end{align}
The terms $\alpha_i^\pm \geq 0$ are the stoichiometric coefficients of the
reactions. Each reaction (precipitation and dissolution) is characterized by a reaction constant $\lambda^\pm$, being
$\lambda^+$ the one associated with the precipitation and  $\lambda^-$ the one
associated with the dissolution. We have $\lambda^\pm \geq 0$. We indicate with $\{u_i\}_{i=1}^{N}$ and $w$ the molar concentration of
the species $\{U_i\}_{i=1}^{N}$ and $W$, respectively. We have
the lower bound $u_i\geq0$, for all $i=1,\ldots,N$, as well as $w\geq0$.
We can write the net precipitation rate associated with
the reaction \eqref{eq:chemical} in the following
way
\begin{gather*}
    r_w(\{u_i\}_{i=1}^{N}) = \lambda^+ \prod_{i=1}^N u_i^{\alpha_i^+} - \lambda^-
    \prod_{i=1}^N u_i^{\alpha_i^-},
\end{gather*}
the first term being the rate of creation of solid part $w$ and the
second term the dissolution rate of the solid part in a unit time.

For simplicity, we suppose only one spice of positive ions and one spice of
negative ions, meaning $N=2$. Moreover, we assume electrical equilibrium, \textit{i.e.} number of anions equal to the
number of cations, and thus we can have $u_i = u$ for $i=1, 2$. The
previously introduced reaction rate can be simplified as
\begin{gather}\label{eq:chemicalsimpl}
    r_w(u) = \lambda^+ u^{\alpha^+} - \lambda^- u^{\alpha^-}
    \quad \text{with} \quad
    \alpha^\pm = \alpha_1^\pm + \alpha_2^\pm.
\end{gather}
We consider that the dissolution of $w$ does not depend on the presence of
ions $u$ whereas precipitation involves all the ions
present. In formula, if we assume
that $\alpha^-_i = 0$, for $i=1, 2$, and \eqref{eq:chemicalsimpl}
becomes
\begin{gather*}%\label{eq:chemicalsimpl}
    r_w(u) = \lambda^-\left( \frac{\lambda^+}{\lambda^-} u^{\alpha^+} -
    1\right).
\end{gather*}
Inspired by the previous relation, we can finally write the more
abstract reaction rate law that is considered in this work. In a more compact way, we have
\begin{gather*}
    r_w(u) = \lambda [r(u) - 1],
\end{gather*}
with $\lambda\geq 0$ a coefficient {and $r(u)=u^\zeta$}, with $\zeta$ a positive
integer.
The previous models suffer of an inconsistency, in fact they do not not vanish in
the limit case of $w=0$ and might create negative values of the quantities
involved. To overcome this problem, {we reformulate the reaction rate as follows}
\begin{gather}\label{eq:reactionn}
    \begin{aligned}
        r_w(u, w) =
        \begin{dcases*}
            \lambda [r(u) - 1] & if $r(u) -1 \geq 0$\\
            - \lambda[r(u) - 1] & if $r(u) -1 < 0$ and $w > 0$\\
            0 & if $r(u) - 1 < 0$ and $w \leq 0$
        \end{dcases*}
    \end{aligned}.
\end{gather}
The first condition models the case {of a positive net precipitation rate, i.e.}
ions precipitate and form the salt $w$. The second condition requires that the
precipitate $w$ is present, i.e.  $w>0$, and allows for its dissolution. The
last equation stops the reaction  when the dissolution should occur
but the precipitate is not present.

The chemical model we are considering is rather general and it does not depend
on the fact that the solutes are transported in a porous media. However, the
case of reactive transport in a porous medium has the unique feature that the
porous medium itself is influenced by the fact that these reactions, occurring in
each spatial point of the domain of interest, can alter its porosity and
permeability.

\begin{figure}[tbp]
    \centering
    \includegraphics[width=0.5\textwidth]{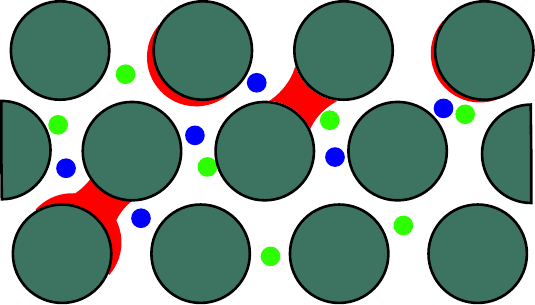}
    \caption{Graphical representation of a porous media in presence of reactive
    species. The floating green and blue circles represent the anions and
    cations flowing in the void space between solid rock grains. The red parts are the
    deposited material due to the reaction.}
    \label{fig:depisition}
\end{figure}

We assume that the solid matrix is formed by two distinct parts: the
precipitate $w$ and the solid inert part that does not react. The latter will be called
solid rock. In the absence
of precipitate the porous media has a prescribed or reference values of porosity and
permeability due to the solid rock, named $\overline{\phi}$ and $\overline{k}$ respectively. During the
flow of chemical species in the porous media, transported by the liquid phase,
a reaction may happen and the
deposition of new material is assumed to be around the grains of solid rock or on a
layer of precipitate already deposited. See \cite{Noorden2009} for a more
detailed discussion. A graphical representation is given in Figure
\ref{fig:depisition}, where we can notice that the deposition of new material
alters the flow path in the porous media itself.
We can model the change of porosity of porous media by a law which accounts for the dependence on the precipitate concentration as follows
\begin{gather}\label{eq:porosity}
    \begin{aligned}
        &\partial_t \phi = - \nu(\phi) \partial_t w && t > 0\\
        &\phi(t=0) = \overline{\phi}
    \end{aligned}
\end{gather}
where the precipitate dependent function, which represents the rate of
deposition of the solute around the solid rock grains, has the properties
\begin{gather*}
    \nu \geq 0 \quad \text{and} \quad \phi = 0 \quad \Rightarrow \quad \nu = 0.
\end{gather*}
We notice that when
$\phi=0$ the porous media is occluded and no deposition of new material can
take place. Moreover, \eqref{eq:porosity} allows the porosity to increase in
presence of the dissolution of precipitate, conversely the porosity decreases
when the precipitate is deposited.
Other model can be taken into consideration, but to keep the
presentation simpler we adopt \eqref{eq:porosity} where $\nu(\phi) = \eta
\phi$, with $\eta$ a positive constant.

Finally, also the permeability of the porous media is influenced by the
reaction. In this work, we consider a Kozeny relationship between the porosity
and the permeability $k$, namely
\begin{gather}\label{eq:perm}
    k(\phi) = \overline{k} \frac{\phi^\alpha}{\overline{\phi}^\alpha},
\end{gather}
with $\alpha > 0$ a rock dependent parameter. In this work we chose $\alpha=2$.
More sophisticated models can be found in, e.g., \cite{Hommel2018}.

\subsection{Transport model}

We introduce now the transport model, assuming that the anions and cations are transported in the porous media as passive scalars. Meaning that there is not
a direct influence of the scalar variable $u$ on the given advective field $\bm{q}$.
In addition, we consider a Fick's law to describe the molecular diffusivity of
$u$ in the liquid with a coefficient {(or tensor)} $d$. The model we are considering for the
solute $u$ is given, in its mixed formulation, by
\begin{gather}\label{eq:transport}
    \begin{aligned}
        &\begin{aligned}
            &\bm{\chi} - \bm{q} u + \phi d \nabla u = \bm{0}\\
            &\partial_t(\phi u ) + \nabla \cdot \bm{\chi} + \phi r_w(u, w) = 0
        \end{aligned}
        && \text{in } \Omega \times \{t >0\}\\
        &\tr u = \hat{u} && \text{on } \Gamma_{in} \times \{t > 0\}\\
        & \tr \phi d \nabla u \cdot \bm{n} = 0 && \text{on } \Gamma_{out}
        \times \{t > 0\}\\
        &\tr \bm{\chi} \cdot \bm{n}= 0 && \text{on } \Gamma_{N} \times \{t > 0\}\\
        &u(t=0) = \overline{u} && \text{in } \Omega
    \end{aligned}.
\end{gather}
{With $\bm{\chi}$ we have denoted the total flux given by the contributions of advection and diffusion.}
The boundary $\partial \Omega$ of the porous media is divided into three
disjoint parts $\Gamma_{in}$, $\Gamma_{out}$, and $\Gamma_{N}$ such that
$\overline{\Gamma_{in}} \cup \overline{\Gamma_{out}} \cup \overline{\Gamma_N} =
\overline{\partial \Omega}$. The portion $\Gamma_{in}$ represents the inflow
boundary, with $\tr \bm{q} \cdot \bm{n} < 0$, where
the value of $u$ is prescribed as $\hat{u}$. The part $\Gamma_{out}$ is where the outflow takes
place with $\tr \bm{q} \cdot \bm{n} > 0$.
On $\Gamma_N$ we prescribe zero flux exchange with the outside, we are
here assuming that $\tr \bm{q} \cdot \bm{n} =0$ {in agreement with the boundary
conditions of the Darcy problem, see Section \ref{subsec:darcy_model}}. Other type of boundary conditions can be considered. The
outward unit normal of $\partial \Omega$ is indicated with $\bm{n}$, and  the operator $\tr$ indicates, in a formal way, a
spacial trace operator mapping the variable at the corresponding portion of the
boundary $\partial \Omega$. Finally, $\overline{u}$ represents the initial data for the solute.

Problem \eqref{eq:transport} is an advection-diffusion-reaction equation, which
can degenerate due to clogging, i.e. $\phi = 0$ from \eqref{eq:porosity}, in some parts
of the domain. The reaction term, described by the law \eqref{eq:reactionn}, is a
non-linear and non-smooth function of the solution $u$ and of the precipitate $w$.

The evolution of the precipitate $w$ follows a similar model of $u$, with the
additional assumption that $w$ does not move in space. All the spatial
differential operators are thus removed and we obtain that
the model  is an ordinary differential equation in each point
of $\Omega$, namely
\begin{gather}\label{eq:precipitate}
    \begin{aligned}
        &\partial_t (\phi w) - \phi r_w(u, w) = 0 && \text{in } \Omega \times
        \{t > 0\}\\
        &w(t=0) = \overline{w} && \text{in }  \Omega
    \end{aligned}.
\end{gather}
The value $\overline{w}$ represents the initial condition of $w$ in $\Omega$.
The reaction terms in the two equations \eqref{eq:transport} and
\eqref{eq:precipitate} {match} each other.

\subsection{Darcy model}\label{subsec:darcy_model}

In this part we introduce the Darcy model and its relation with the previously discussed chemical
model. We are interested in computing the Darcy velocity
$\bm{q}$ and the pressure field $p$ in the porous media satisfying the following
relations
\begin{gather}\label{eq:darcy}
    \begin{aligned}
        &\begin{aligned}
            & \bm{q} = - k(\phi) \nabla p \\
            &\partial_t \phi + \nabla \cdot \bm{q} = f
        \end{aligned}
        && \text{in } \Omega \times \{t>0\}\\
        &\tr p = \overline{p} && \Gamma_{in} \times \{t>0\}\\
        &\tr \bm{q}\cdot\bm{n} = \overline{q} && \Gamma_{out} \times \{t>0\}\\
        &\tr \bm{q} \cdot \bm{n} = 0 && \Gamma_N \times \{t>0\}
    \end{aligned}.
\end{gather}
The division of the boundary $\partial \Omega$ into parts follows the
description given in \eqref{eq:transport}.
The boundary value $\overline{p}$ represents the data at the inflow. The value
$\overline{q}$ is the outflow flux out of $\Gamma_{out}$ with the request that
$\overline{q}>0$. The condition on $\Gamma_N$ is a no flow condition for that
portion of boundary.
By conservation we obtain that $\tr \bm{q}\cdot \bm{n} < 0$ on $\Gamma_{in}$.
Also in this case other type of boundary conditions can be considered,  but they should
be coherent with the one prescribed in the model \eqref{eq:transport}.

Equation \eqref{eq:darcy} is coupled with the reactive models
\eqref{eq:transport} and \eqref{eq:precipitate} via the dependency of porosity
and permeability on the solute $u$ and precipitate $w$.

\subsection{The complete model}\label{subsec:complete}

The complete model is a six {unknowns} model and describes the evolution in time and space of:
\textit{i)} $u$ solute, \textit{ii)} $w$ precipitate, \textit{iii)} $\phi$
porosity, \textit{iv)} $k$
permeability, \textit{v)} $\bm{q}$ Darcy velocity, and \textit{vi)} $p$ pressure. The equations
involved are \eqref{eq:transport}, \eqref{eq:precipitate}, \eqref{eq:porosity},
\eqref{eq:perm}, and \eqref{eq:darcy} respectively for each variable or pair of
variables. The
resulting system is fully coupled, non-smooth and non-linear with possible
degeneracy due to vanishing porosity and permeability.

\section{A reduced model a fracture}\label{sec:reducedmodel}

A fracture is a thin object immersed in a porous media, whose aperture is orders
of magnitude smaller than any other characteristic size of the problem at hand.
The fact that the fracture may exhibit higher or lower permeability with respect to the
surrounding porous media increases the problem complexity and requires a proper
treatment to obtain an effective and reliable model. The choice adopted here is
a reduced model, meaning that the fracture is reduced as an object of lower dimension and new equations and coupling conditions are derived.

In this part, we start by presenting the interface conditions used to couple the
porous media and the fracture, being the latter represented as an equi-dimensional object.
Then, we present the model reduction procedure to introduce the new model and
interface conditions.

\subsection{Coupling condition for the equi-dimensional model}

\begin{figure}[tbp]
    \centering
    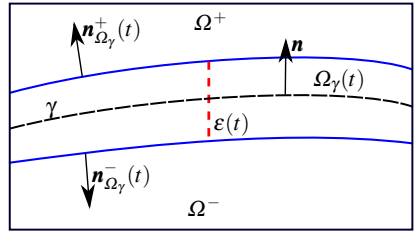
    \caption{Equi-dimensional representation of a fracture $\Omega_\gamma$ immersed in a porous
    media $\Omega$.}
    \label{fig:fracture}
\end{figure}

Given a parameter $\epsilon(t)$, called the fracture aperture, which might
change in time due to deposition of dissolution of new material. Following the
presentation given in \cite{Formaggia2012} we can define the
fracture as the domain $\Omega_\gamma(t)$ given by
\begin{gather}\label{eq:fracture}
    \Omega_\gamma(t) = \left\{ \bm{x} \in \mathbb{R}^n: \, \bm{x} = \bm{s} +
    \xi(t)
    \bm{n}, \text{ with } \bm{s} \in \gamma \text{ and } \xi \in \left(
    -\frac{\epsilon(\bm{s}, t)}{2}, \frac{\epsilon(\bm{s}, t)}{2}
    \right)\right\},
\end{gather}
where $\gamma$ is a non self-intersecting one-codimensional manifold of class
$C^2$. We have $\epsilon(\cdot, t) \in C^2(\gamma)$ and we assume that the fracture
aperture varies slowly compare to the local coordinate system. The vector $\bm{n}$ is the
normal vector of $\gamma$ pointing towards one of the side of the surrounding porous
media. This choice of orientation is arbitrary and will not change the following procedure. Finally, to
ease the presentation we suppose that the fracture cuts the porous media in two
disjoint parts indicated with $+$ and $-$. Extension to more general cases are
straightforward.
An example is reported in Figure \ref{fig:fracture}.

Being $\Omega_\gamma$ equi-dimensional with respect to the surrounding porous
media $\Omega$ it is possible to write the same equations to model the reactive
transport as the one discussed in Subsection \ref{subsec:complete} but applied
to $\Omega_\gamma$ instead. We will indicate with a subscript if the variable or
data is
referred to the porous media $\Omega$ or to the equi-dimensional representation
of the fracture $\Omega_\gamma$.

In addition to this, interface conditions have to be
considered to couple the two problems at their common boundaries. For the
transport equation \eqref{eq:transport}, following \cite{Fumagalli2012a}, we have the conservation of the total
flux and the continuity of the solute $u$, meaning
\begin{align}\label{eq:coupling_transport}
    &
    \begin{aligned}
        &\tr \bm{\chi}_\Omega \cdot \bm{n}_{\Omega_\gamma} = \tr
        \bm{\chi}_{\Omega_\gamma} \cdot \bm{n}_{\Omega_\gamma}\\
        &\tr u_\Omega = \tr u_{\Omega_\gamma}
    \end{aligned}
    && \text{ on } \partial \Omega \cap \partial \Omega_\gamma,
\end{align}
where $\bm{n}_{\Omega_\gamma}$ is the unit normal of the boundary of
$\Omega_\gamma$ pointing from the latter toward $\Omega$.
For the Darcy equation \eqref{eq:darcy}  across the
interfaces we have continuity of the normal component of Darcy velocity $\bm{q}$ as well as the continuity of the
pressure $p$. Following \cite{Martin2005,Faille2014a,Schwenck2015} we obtain
\begin{align}\label{eq:coupling_darcy}
    &
    \begin{aligned}
        &\tr \bm{q}_\Omega \cdot \bm{n}_{\Omega_\gamma} = \tr
        \bm{q}_{\Omega_\gamma} \cdot \bm{n}_{\Omega_\gamma}\\
        &\tr p_\Omega = \tr p_{\Omega_\gamma}
    \end{aligned}
    && \text{ on } \partial \Omega \cap \partial \Omega_\gamma.
\end{align}

Finally, the full equi-dimensional model for porous media-fracture system is
given by equations \eqref{eq:transport}, \eqref{eq:precipitate}, \eqref{eq:porosity},
\eqref{eq:perm}, and \eqref{eq:darcy} for both $\Omega$ and $\Omega_\gamma$
along with the coupling conditions given by \eqref{eq:coupling_transport} and
\eqref{eq:coupling_darcy}.

\subsection{The reduced variables}

\begin{figure}[tbp]
    \centering
    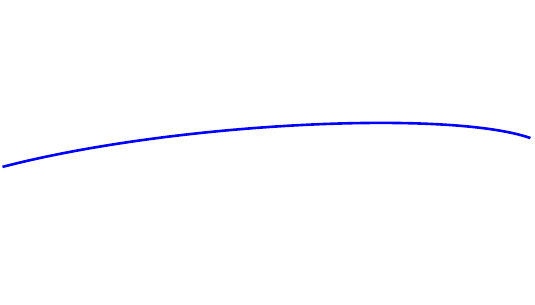
    \caption{Hybrid-dimensional representation of a fracture immersed in a porous
    media.}
    \label{fig:fracture_reduced}
\end{figure}

The model reduction procedure approximates the equi-dimensional representation of
the fracture $\Omega_\gamma$ by its centre line $\gamma$, and derive new
equations to describe the variables in $\gamma$ and new interface conditions for
the
coupling with the surrounding porous media. Due to the previously mentioned
assumptions on $\gamma$, we approximate $\bm{n}_{\Omega_\gamma}^\pm$ with $\pm
\bm{n}$. The representation of $\Omega \subset \mathbb{R}^n$ and $\gamma$ as
co-dimension one object is usually named as hybrid-dimensional. See Figure \ref{fig:fracture_reduced} as an example.

The new variables defined on $\gamma$ are defined differently if they are scalar
or vector fields. In the former case we {define average values as}
\begin{gather}\label{eq:reduced_variables1}
    u_\gamma(\bm{s}, t) = \frac{1}{\epsilon(\bm{s}, t)}
    \int_{-\frac{\epsilon(\bm{s}, t)}{2}}^{\frac{\epsilon(\bm{s}, t)}{2}}
    u_{\Omega_\gamma}(t)
    d \bm{n}(\bm{s})
    \quad \text{and} \quad
    p_\gamma(\bm{s}) = \frac{1}{\epsilon(\bm{s}, t)}
    \int_{-\frac{\epsilon(\bm{s}, t)}{2}}^{\frac{\epsilon(\bm{s}, t)}{2}}
    p_{\Omega_\gamma}(t)
    d \bm{n}(\bm{s}),
\end{gather}
where $\bm{s} \in \gamma$ and the time dependent integrals are done along the  direction normal to the fracture $\gamma$. For the vector fields $\bm{\chi}$ and $\bm{q}$ we need
to introduce the following projection matrices along and across the fracture,
given by
\begin{gather*}
    N = \bm{n} \otimes \bm{n}
    \quad \text{and} \quad
    T = I - N.
\end{gather*}
Now, we can define
\begin{gather}\label{eq:reduced_variables2}
    \bm{\chi}_\gamma =
    \int_{-\frac{\epsilon(\bm{s}, t)}{2}}^{\frac{\epsilon(\bm{s}, t)}{2}}
    T(\bm{s}) \bm{\chi}_{\Omega_\gamma}(t) d \bm{n}(\bm{s})
    \quad \text{and} \quad
    \bm{q}_\gamma =
    \int_{-\frac{\epsilon(\bm{s}, t)}{2}}^{\frac{\epsilon(\bm{s}, t)}{2}}
    T(\bm{s}) \bm{q}_{\Omega_\gamma}(t) d \bm{n}(\bm{s})
\end{gather}
which, it is worth to notice, are not average values of fluxes or velocity, but
time dependent integrals.
We require that, following the idea of \cite{Alboin2002,Martin2005}, the
permeability for the fracture  \eqref{eq:darcy} is aligned to the local
coordinate system,
meaning that we have
\begin{gather} \label{eq:assump_perm}
    k_{\Omega_\gamma} = k_\gamma T + \kappa N,
\end{gather}
where $k_\gamma$ is a square tensor of dimension $n-1$, symmetric and positive defined,
and $\kappa$ is a positive real number. Moreover, also the diffusion tensor
of equation \eqref{eq:transport} is required to follow the similar request
\begin{gather*}
    d_{\Omega_\gamma} = d_\gamma T + \delta N,
\end{gather*}
where $d_\gamma$ is a square tensor of dimension $n-1$, symmetric and positive
defined, and $\delta$ is a positive real number.

Finally, the fracture is initially considered open, meaning $\phi_{\Omega_\gamma}
= 1$, and in the reduced model its role will be played by the aperture
$\epsilon$. We will give a specific law for its evolution. This will be part of the discussion in
Subsection \ref{subsub:aperture}.

\subsection{Reduced transport model}\label{subsub:reduced_transport}

We describe now the procedure to derive the reduced model for the system
\eqref{eq:transport}. First of all the first equation of \eqref{eq:transport} is decomposed in its
normal and tangential parts as
\begin{gather}\label{eq:transport_0}
    \begin{aligned}
        &T \bm{\chi}_{\Omega_\gamma} - T \bm{q}_{\Omega_\gamma} u_{\Omega_\gamma} +  T d_{\Omega_\gamma}
        \nabla u_{\Omega_\gamma} = \bm{0},\\
        &N \bm{\chi}_{\Omega_\gamma} - N \bm{q}_{\Omega_\gamma}
        u_{\Omega_\gamma} +  N d_{\Omega_\gamma}
        \nabla u_{\Omega_\gamma} = \bm{0}.
    \end{aligned}
\end{gather}
Now, the tangential equation is integrated in the normal direction $\bm{n}$ across the
fracture. Dropping the dependency on $\bm{s}$ and $t$ when not needed, we obtain
\begin{gather*}
    \int_{-\frac{\epsilon}{2}}^{\frac{\epsilon}{2}}
    T \bm{\chi}_{\Omega_\gamma} d \bm{n}
    -
    \int_{-\frac{\epsilon}{2}}^{\frac{\epsilon}{2}}
    T \bm{q}_{\Omega_\gamma} u_{\Omega_\gamma} d \bm{n}
    +
    \int_{-\frac{\epsilon}{2}}^{\frac{\epsilon}{2}}
    T d_{\Omega_\gamma} \nabla u_{\Omega_\gamma} d
    \bm{n}= \bm{0},
\end{gather*}
having $T d_{\Omega_\gamma} = T d_\gamma$ and by assuming small variations of along the thickness of the fracture of
$\bm{q}_\gamma$ and $u_\gamma$,
we get the following expression
\begin{gather}\label{eq:transport_reduced_1}
    \bm{\chi}_\gamma
    -
    \bm{q}_\gamma u_\gamma
    +
    \epsilon d_\gamma \nabla_T u_\gamma = \bm{0} \qquad \text{in } \gamma \times
        \{t>0\},
\end{gather}
where the nabla operator $\nabla_T = T \nabla $ is defined now on the tangent space of the fracture.
The second equation of \eqref{eq:transport_0} gives the coupling conditions
between the fracture $\gamma$ and the surrounding porous medium, i.e. the sides
$\Omega^+$ and $\Omega^-$. The derivation of such conditions requires the
integration from the centre line of $\Omega_\gamma$ to its boundary, which is
given, for $\Omega^+$, by
\begin{gather*}
    \int_{0}^{\frac{\epsilon}{2}}
    N \bm{\chi}_{\Omega_\gamma} \cdot \bm{n} d \bm{n}
    -
    \int_{0}^{\frac{\epsilon}{2}}
    N \bm{q}_{\Omega_\gamma} u_{\Omega_\gamma} \cdot \bm{n} d \bm{n}
    +
    \int_{0}^{\frac{\epsilon}{2}}
    N d_{\Omega_\gamma} \nabla u_{\Omega_\gamma} \cdot \bm{n} d
    \bm{n}= 0.
\end{gather*}
For the last term, we have $N d_{\Omega_\gamma} = N \delta$ and it can be approximated as
\begin{gather*}
    \int_{0}^{\frac{\epsilon}{2}}
    N \delta \nabla u_{\Omega_\gamma} \cdot \bm{n} d
    \bm{n} \approx \delta (\tr u_{\Omega^+} - u_\gamma).
\end{gather*}
While for the other two terms we consider a first order one-side integration
rule, along with the continuity conditions \eqref{eq:coupling_transport}  and
\eqref{eq:coupling_darcy} for its approximation. We get
\begin{gather*}
    \int_{0}^{\frac{\epsilon}{2}}
    N \bm{\chi}_{\Omega_\gamma} \cdot \bm{n} d \bm{n} -
    \int_{0}^{\frac{\epsilon}{2}}
    N \bm{q}_{\Omega_\gamma} u_{\Omega_\gamma} \cdot \bm{n} d \bm{n}\approx
    \frac{\epsilon}{2}\left(
    \tr \bm{\chi}_{\Omega^+} \cdot \bm{n}
    -
    \tr \bm{q}_{\Omega^+} \cdot \bm{n} \tr u_{\Omega^+}\right).
\end{gather*}
Finally, we get the coupling condition for the side of the fracture in contact
with $\Omega^+$
\begin{gather}\label{eq:transport_coupling_reduced}
    \epsilon \left( \tr \bm{\chi}_{\Omega^+} \cdot \bm{n}
    -
    \tr \bm{q}_{\Omega^+} \cdot \bm{n} \tr u_{\Omega^+} \right) =
    2 \delta (u_\gamma - \tr u_{\Omega^+})
\end{gather}
For the other side $\Omega^-$ the derivation is similar.

\begin{figure}[tbp]
    \centering
    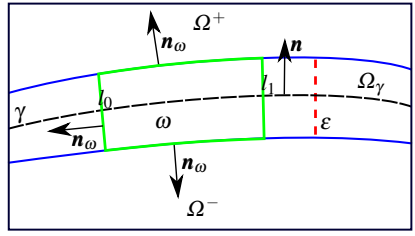
    \caption{Equi-dimensional representation of a fracture immersed in a porous
    media with the control volume $\omega$.}
    \label{fig:fracture_omega}
\end{figure}

The conservation equation, second of \eqref{eq:transport},  is reduced following
the same approach presented in \cite{Boon2018}. Its  integral
form is given by
\begin{gather}\label{eq:cons}
    \partial_t
    \int_{\omega(t)}
     u_{\Omega_\gamma} d \bm{x}
    +
    \int_{\partial \omega(t)}
    \tr \bm{\chi}_{\Omega_\gamma} \cdot \bm{n}_{\omega} d \bm{\sigma}
    +
    \int_{\omega(t)}
    r_w(u_{\Omega_\gamma}, w_{\Omega_\gamma})  d \bm{x}= 0
\end{gather}
where $\omega(t) = (l_0, l_1) \times (-{\epsilon(t)}/{2},
{\epsilon(t)}/{2}) \subset\Omega_\gamma(t)$ and with $(l_0, l_1) \subset
\gamma$.
Note that the latter does not depend on time. The vector $\bm{n}_\omega$ is the outward unit
normal of $\omega$. See
Figure \ref{fig:fracture_omega} for a more detailed representation of the
objects involved.
The first and third part of the previous equation, omitting
the dependency on $t$, are now
given by
\begin{gather*}
    \int_{l_0}^{l_1} \left(
    \partial_t
    \int_{-\frac{\epsilon}{2}}^{\frac{\epsilon}{2}}
     u_{\Omega_\gamma} d \bm{n}
    +
    \int_{-\frac{\epsilon}{2}}^{\frac{\epsilon}{2}}
    r_w(u_{\Omega_\gamma}, w_{\Omega_\gamma})  d \bm{n} \right) d \bm{s}=
    \int_{l_0}^{l_1} \partial_t( \epsilon u_\gamma) + \epsilon r_w(u_{\gamma}, w_{\gamma})d \bm{s}
\end{gather*}
The second part of \eqref{eq:cons} becomes
\begin{gather*}
    \int_{\partial \omega} \tr \bm{\chi}_{\Omega_\gamma} \cdot
    \bm{n}_{\omega} d \bm{\sigma}  = \int_{\partial \Omega_\gamma \cap
    \partial \omega}
    \tr \bm{\chi}_{\Omega_\gamma} \cdot
    \bm{n}_{\Omega_\gamma} d \bm{\sigma} + \int_{\partial \omega^+}
    \tr \bm{\chi}_{\Omega_\gamma} \cdot
    \bm{n}_{\omega} d \bm{\sigma}
    +\\ \int_{\partial \omega^-}
    \tr \bm{\chi}_{\Omega_\gamma} \cdot
    \bm{n}_{\omega} d \bm{\sigma},
\end{gather*}
with $\partial \omega^+ = \{l_1\}\times(-\epsilon/2,
\epsilon/2) $ and $\partial \omega^- = \{l_0\}\times(-\epsilon/2,
\epsilon/2) $. Now, setting $\abs{(l_0, l_1)} \rightarrow 0$ we obtain the following expressions
\begin{gather*}
    \lim_{\abs{(l_0, l_1)} \rightarrow 0} \frac{1}{\abs{(l_0, l_1)}} \int_{l_0}^{l_1} \partial_t( \epsilon u_\gamma) +
    \epsilon r_w(u_{\gamma}, w_{\gamma})d \bm{s} = \partial_t( \epsilon u_\gamma) +
    \epsilon r_w(u_{\gamma}, w_{\gamma})\\
    \lim_{\abs{(l_0, l_1)} \rightarrow 0}\frac{1}{\abs{(l_0, l_1)}}\int_{\partial \omega} \tr \bm{\chi}_{\Omega_\gamma} \cdot
    \bm{n}_{\omega} d \bm{\sigma} = \tr \bm{\chi}_{\Omega_\gamma} \cdot
    \bm{n}_{\Omega_\gamma}|_{\frac{\epsilon}2} +
    \tr \bm{\chi}_{\Omega_\gamma} \cdot
    \bm{n}_{\Omega_\gamma}|_{-\frac{\epsilon}2}
    + \nabla_T \cdot \bm{\chi}_\gamma
\end{gather*}
by using the continuity conditions \eqref{eq:coupling_transport} the last
equation becomes
\begin{gather*}
    \lim_{\abs{(l_0, l_1)} \rightarrow 0} \frac{1}{\abs{(l_0, l_1)}}\int_{\partial \omega} \tr \bm{\chi}_{\Omega_\gamma} \cdot
    \bm{n}_{\omega} d \bm{\sigma} = \tr \bm{\chi}_{\Omega^+} \cdot
    \bm{n} -
    \tr \bm{\chi}_{\Omega^-} \cdot
    \bm{n}
    + \nabla_T \cdot \bm{\chi}_\gamma.
\end{gather*}
Finally, the conservation equation for the transport system of the solute is reduced as
\begin{gather}\label{eq:reduced_solute}
    \begin{aligned}
        &\partial_t( \epsilon u_\gamma)    + \nabla_T \cdot \bm{\chi}_\gamma +
        \tr \bm{\chi}_{\Omega^+} \cdot
        \bm{n} -
        \tr \bm{\chi}_{\Omega^-} \cdot
        \bm{n} +
        \epsilon r_w(u_{\gamma}, w_{\gamma}) = 0
        && \text{in } \gamma \times
        \{t>0\}\\
        &u_\gamma(t=0) = \overline{u}_\gamma && \text{in } \gamma
    \end{aligned},
\end{gather}
where $\overline{u}_\gamma$ is the reduced initial condition, given by
$\overline{u}_\gamma = \frac{1}{\epsilon}
\int_{-\frac{\epsilon}{2}}^{\frac{\epsilon}{2}} \overline{u}_{\Omega_\gamma}$.
The reduced boundary conditions for the transport problem are given by the following
\begin{gather}\label{eq:reduced_solute_bc}
    \begin{aligned}
        &\tr u_\gamma = \hat{u}_\gamma && \text{on } (\partial \gamma \cap \Gamma_{in}) \times \{t > 0\}\\
        & \tr \epsilon d_\gamma \nabla_T u_\gamma \cdot \bm{n} = 0 && \text{on } (\partial \gamma
        \cap \Gamma_{out})
        \times \{t > 0\}\\
        &\tr \bm{\chi}_\gamma \cdot \bm{n}= 0 && \text{on } (\partial \gamma \cap
        \Gamma_{N}) \times \{t > 0\}
    \end{aligned},
\end{gather}
where $\hat{u}_\gamma$ is defined accordingly. We have assumed here that, for
example if $\gamma$ is one-dimensional, a single boundary condition is assigned
to each end point $\partial \gamma$.

For the precipitate the derivation of the reduced model is rather easy since no
spatial differential operators are involved. From \eqref{eq:precipitate} and by
considering again the control volume $\omega(t)$, we get
\begin{gather*}
    \partial_t\int_{\omega(t)}
    w_{\Omega_\gamma}  d \bm{x}-
    \int_{\omega(t)}
    r_w(u_{\Omega_\gamma}, w_{\Omega_\gamma})  d \bm{x}= 0
\end{gather*}
and proceeding as before we obtain
\begin{gather*}
    \lim_{\abs{(l_0, l_1)} \rightarrow 0} \frac{1}{\abs{(l_0, l_1)}}\int_{l_0}^{l_1} \left(
    \partial_t
    \int_{-\frac{\epsilon}{2}}^{\frac{\epsilon}{2}}
    w_{\Omega_\gamma} d \bm{n}
    -
    \int_{-\frac{\epsilon}{2}}^{\frac{\epsilon}{2}}
    r_w(u_{\Omega_\gamma}, w_{\Omega_\gamma})  d \bm{n} \right) d \bm{s}=\\
    \lim_{\abs{(l_0, l_1)} \rightarrow 0} \frac{1}{\abs{(l_0, l_1)}}
    \int_{l_0}^{l_1} \partial_t( \epsilon w_\gamma) - \epsilon r_w(u_{\gamma},
    w_{\gamma})d \bm{s} =\partial_t (\epsilon w_\gamma) - \epsilon r_w(u_\gamma,
    w_\gamma).
\end{gather*}
Finally, the following is the reduced model for the precipitate $w$
\begin{gather}\label{eq:reduced_precipitate}
    \begin{aligned}
        &\partial_t (\epsilon w_\gamma) - \epsilon r_w(u_\gamma, w_\gamma) = 0
        && \text{in } \gamma \times \{t>0\}\\
        &w_\gamma(t=0) = \overline{w}_\gamma&& \text{in } \gamma
    \end{aligned},
\end{gather}
where $\overline{w}_\gamma$ is the reduced initial condition for $w_\gamma$,
given by $\overline{w}_\gamma = \frac{1}{\epsilon}
\int_{-\frac{\epsilon}{2}}^{\frac{\epsilon}{2}} \overline{w}_{\Omega_\gamma}$.

\subsection{Aperture and permeability models}\label{subsub:aperture}

\begin{figure}[tbp]
    \centering
    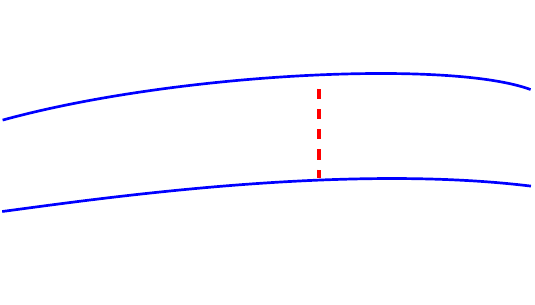
    \caption{Equi-dimensional representation of a fracture immersed in a porous
    media with dynamics of deposition and dissolution due to the chemical
    reaction. The solutes are the blue and green circles and the precipitate is
    depicted in red.}
    \label{fig:fracture_aperture}
\end{figure}

Following the ideas discussed in \cite{Kumar2011},
to derive the variation of the fracture aperture by the deposition or
dissolution of the solute, we consider a law similar to the one given for the
porosity in \eqref{eq:porosity}. However, in this case since the fracture is
supposed to be initially empty, free from granular material, we assume that the new
material is accumulated or dissolved at the fracture boundary. See Figure
\ref{fig:fracture_aperture} for a graphical representation.
We consider again a precipitate dependent law to describe the rate of aperture
change, we have
\begin{gather}\label{eq:aperture}
    \begin{aligned}
        &\partial_t \epsilon = - \upsilon(\epsilon) \partial_t w_\gamma && t > 0\\
        &\epsilon(t=0) = \overline{\epsilon}
    \end{aligned}
\end{gather}
where the aperture dependent model, which represents the rate of deposition of
the solute around the fracture walls, has the following properties
\begin{gather*}
    \upsilon \geq 0 \quad \text{and} \quad \epsilon = 0 \quad \Rightarrow \quad
    \upsilon = 0.
\end{gather*}
We notice that when
$\epsilon=0$ the fracture is occluded and no deposition of new material
takes place.
Moreover, \eqref{eq:aperture} allows the aperture to increase in
presence of the dissolution of precipitate, conversely the aperture decreases
when the precipitate is deposited.
Other models can be taken into consideration, but to keep the
presentation simpler we adopt \eqref{eq:aperture} where $\upsilon(\epsilon) =
\eta_\gamma
\epsilon$, with $\eta_\gamma$ a positive constant. In \eqref{eq:aperture}, the
value of $\overline{\epsilon} \geq 0$ represents the initial aperture of the
fracture.

The fracture permeability, both normal $\kappa$ and tangential $k_\gamma$ are
now related to the fracture aperture by the cubic law
\begin{gather}\label{eq:reduced_permeability}
    k_\gamma = \overline{k}_\gamma \frac{\epsilon^2}{\overline{\epsilon}}
    \quad \text{and} \quad
    \kappa = \overline{\kappa} \frac{\epsilon^2}{\overline{\epsilon}}.
\end{gather}
Here $\overline{k}_\gamma$, symmetric and positive defined, and
$\overline{\kappa} > 0 $ are the reference tangential and normal fracture
permeability, respectively.

\subsubsection{Reduced Darcy model}\label{subsub:reduced_darcy}

In this part we derive the reduced model for the Darcy system \eqref{eq:darcy}
written in the fracture. The steps are rather similar to the one presented for
the transport equation with few modifications. The Darcy equation, first of
\eqref{eq:darcy} is projected on the tangential and normal directions of the
fracture obtaining
\begin{gather}\label{eq:darcy_normal_tange}
    \begin{aligned}
        T \bm{q}_{\Omega_\gamma} + T k_{\Omega_\gamma} \nabla p_{\Omega_\gamma} = \bm{0},\\
        N \bm{q}_{\Omega_\gamma} + N k_{\Omega_\gamma} \nabla p_{\Omega_\gamma} = \bm{0}.
    \end{aligned}
\end{gather}
The first of \eqref{eq:darcy_normal_tange} is now integrated across the normal section of the fracture,
along the direction given by $\bm{n}$. We have
\begin{gather*}
    \int_{-\frac{\epsilon}{2}}^{\frac{\epsilon}{2}}
    T \bm{q}_{\Omega_\gamma} d \bm{n}
    +
    \int_{-\frac{\epsilon}{2}}^{\frac{\epsilon}{2}}
    T k_{\Omega_\gamma} \nabla p_{\Omega_\gamma} d \bm{n}= \bm{0}.
\end{gather*}
From the assumption on the permeability \eqref{eq:assump_perm} we obtain $T
k_{\Omega_\gamma} = T k_\gamma$. By assuming small variations along the
thickness of the fracture of $\nabla p_{\Omega_\gamma}$, we get the following
relation
\begin{gather}\label{eq:darcy_reduceddd}
    \bm{q}_\gamma + \epsilon k_\gamma \nabla_T p_\gamma = \bm{0}
    \quad \text{in } \quad
    \gamma \times \{t>0\}.
\end{gather}
The second relation in \eqref{eq:darcy_normal_tange} gives the coupling
conditions between the fracture and the surrounding porous media for the Darcy
problem. The approach is similar to the one already presented for the transport
part, we integrate the second of \eqref{eq:darcy_normal_tange} from $0$ to
$\epsilon/2$ and we do some approximation of the integrals involved. We start
with
\begin{gather*}
    \int_{0}^{\frac{\epsilon}{2}}
    N \bm{q}_{\Omega_\gamma} \cdot \bm{n} +
    \int_{0}^{\frac{\epsilon}{2}}
    N k_{\Omega_\gamma} \nabla p_{\Omega_\gamma} \cdot \bm{n}
    = \bm{0}.
\end{gather*}
The first integral is approximated by a one-side integration rule and the
coupling conditions \eqref{eq:coupling_darcy} are considered to get
\begin{gather*}
    \int_{0}^{\frac{\epsilon}{2}}
    N \bm{q}_{\Omega_\gamma} \cdot \bm{n} \approx \frac{\epsilon}{2} \tr
    \bm{q}_{\Omega^+} \cdot \bm{n},
\end{gather*}
while recognising that $ N k_{\Omega_\gamma} = N \kappa$, for the second term we obtain
\begin{gather*}
    \int_{0}^{\frac{\epsilon}{2}}
    N k_{\Omega_\gamma} \nabla p_{\Omega_\gamma} \cdot \bm{n}\approx
    \kappa (\tr p_{\Omega^+} - p_\gamma).
\end{gather*}
The coupling conditions of the reduced model for the Darcy equation,
for the side of the fracture in contact with
$\Omega^+$, are thus
given by
\begin{gather}\label{eq:darcy_coupling_reduced}
    \epsilon \tr \bm{q}_{\Omega^+} \cdot \bm{n} = 2\kappa (\tr p_{\Omega^+} - p_\gamma).
\end{gather}
Also in this case, for the side $\Omega^-$ the derivation of the coupling
conditions are similar.

Finally, to complete the Darcy system the conservation equation for the fracture
has to be reduced. {Unlike the previous steps, which are in agreement with the
existing literature, see for instance \cite{Martin2005,Formaggia2012,Boon2018}, this last step differ from the
previous works on model reduction for fractured media because we have to account
for a time dependent aperture.} We consider again the control volume $\omega(t)$
given as
before, and the integral form of the conservation equation is given by
\begin{gather}\label{eq:darcy_reducedd}
    \partial_t \int_{\omega(t)} d \bm{x} +
    \int_{\partial \omega(t)} \tr \bm{q}_{\Omega_\gamma} \cdot \bm{n}_\omega d
    \bm{\sigma} = \int_{\omega(t)} f d \bm{x}.
\end{gather}
Reminding that $\omega(t) = (l_0, l_1) \times (-{\epsilon(t)}/{2},
{\epsilon(t)}/{2})$, the first and last terms, dropping the dependence on $t$, can be expressed as
\begin{gather*}
    \lim_{\abs{(l_0, l_1)} \rightarrow 0} \frac{1}{\abs{(l_0, l_1)}}\int_{l_0}^{l_1} \left( \partial_t
    \int_{-\frac{\epsilon}{2}}^{\frac{\epsilon}{2}}
    d \bm{n} - \int_{-\frac{\epsilon}{2}}^{\frac{\epsilon}{2}}
    f d \bm{n}\right) d \bm{s} = \partial_t \epsilon - \epsilon f_\gamma,
\end{gather*}
where $f_\gamma$ is the reduced source or sink term expressed by $f_\gamma =
\frac{1}{\epsilon} \int_{-\frac{\epsilon}{2}}^{\frac{\epsilon}{2}} f d \bm{n}$.
The second term in \eqref{eq:darcy_reducedd} becomes
\begin{gather*}
    \int_{\partial \omega} \tr \bm{q}_{\Omega_\gamma} \cdot \bm{n}_\omega d
    \bm{\sigma}  = \int_{\partial \Omega_\gamma \cap \partial \omega} \tr \bm{q}_{\Omega_\gamma}
    \cdot \bm{n}_\omega d \bm{\sigma} +
    \int_{\partial \omega^+} \tr \bm{q}_{\Omega_\gamma}  \cdot \bm{n}_\omega d
    \bm{\sigma} +\\
    \int_{\partial \omega^-} \tr \bm{q}_{\Omega_\gamma}  \cdot \bm{n}_\omega d
    \bm{\sigma}.
\end{gather*}
By shrinking the domain $\omega$ as $\abs{(l_0, l_1)} \rightarrow 0$ the last
relation becomes
\begin{gather*}
    \lim_{\abs{(l_0, l_1)}\rightarrow 0} \frac{1}{\abs{(l_0, l_1)}}\int_{\partial \omega }
    \tr \bm{q}_{\Omega_\gamma} \cdot \bm{n}_\omega d
    \bm{\sigma} = \tr \bm{q}_{\Omega_\gamma} \cdot
    \bm{n}_\omega|_{\frac{\epsilon}{2}} + \tr \bm{q}_{\Omega_\gamma} \cdot
    \bm{n}_\omega|_{-\frac{\epsilon}{2}} + \nabla_T \cdot \bm{q}_\gamma,
\end{gather*}
and by using the continuity condition at the fracture-porous media boundary
\eqref{eq:coupling_darcy} we
finally get
\begin{gather*}
    \lim_{\abs{(l_0, l_1)}\rightarrow 0} \frac{1}{\abs{(l_0, l_1)}}\int_{\partial \omega }
    \tr \bm{q}_{\Omega_\gamma} \cdot \bm{n}_\omega d
    \bm{\sigma} = \tr \bm{q}_{\Omega^+} \cdot
    \bm{n} - \tr \bm{q}_{\Omega^-} \cdot
    \bm{n}+ \nabla_T \cdot \bm{q}_\gamma.
\end{gather*}
To conclude the conservation equation for the Darcy flow is given by
\begin{gather}\label{eq:reduced_darcy_conservation}
    \partial_t \epsilon + \nabla_T \cdot \bm{q}_\gamma + \tr \bm{q}_{\Omega^+} \cdot
    \bm{n} - \tr \bm{q}_{\Omega^-} \cdot
    \bm{n} = \epsilon f_\gamma \quad
    \text{in } \gamma \times \{t>0\}.
\end{gather}
The reduced boundary conditions for the Darcy problem are given by the following
\begin{gather}\label{eq:reduced_darcy_conservation_bc}
    \begin{aligned}
        &\tr p_\gamma = \overline{p}_\gamma && (\partial \gamma \cap
        \Gamma_{in}) \times \{t>0\}\\
        &\tr \bm{q}_\gamma \cdot\bm{n} = \overline{q}_\gamma && (\partial \gamma
        \cap \Gamma_{out}) \times \{t>0\}\\
        &\tr \bm{q}_\gamma \cdot \bm{n} = 0 && (\partial \gamma \cap \Gamma_N) \times \{t>0\}
    \end{aligned},
\end{gather}
with $\overline{p}_\gamma$ and $\overline{q}_\gamma$ being defined accordingly.

\subsubsection{The complete reduced model}

We can now summarize the full hybrid-dimensional problem, in this case we have six
fields for the porous media and seven other fields for the fractures. For the
former the reader can refer to the description given in Subsection
\ref{subsec:complete} for $\Omega$, which in the latter we have the evolution
of: \textit{i)} $u_\gamma$ solute, \textit{ii)} $w_\gamma$ precipitate,
\textit{iii)} $\epsilon$
aperture, \textit{iv)} $\kappa$ and $k_\gamma$
normal and tangential permeability, \textit{v)} $\bm{q}_\gamma$ Darcy velocity,
and \textit{vi)} $p_\gamma$ pressure.

For the fracture the equations involved are \eqref{eq:transport_reduced_1},
\eqref{eq:transport_coupling_reduced}, \eqref{subsec:complete}, and \eqref{eq:reduced_solute_bc} for
$u_\gamma$. For $w_\gamma$ the problem \eqref{eq:reduced_precipitate}, and for
$\epsilon$ \eqref{eq:aperture}. For the permeabilities $\kappa$ and $k_\gamma$
the model given by \eqref{eq:reduced_permeability}. Finally, for $\bm{q}_\gamma$
and $p_\gamma$ the equations \eqref{eq:darcy_reduceddd}, \eqref{eq:darcy_coupling_reduced},
\eqref{eq:reduced_darcy_conservation}, and
\eqref{eq:reduced_darcy_conservation_bc}.

Finally, it is important to mention that due to the model reduction procedure the
aperture is now a {time dependent} model parameter and not any more a geometrical constraint for the problem.

\section{Conclusion}\label{sec:conclusion}

In this work we have presented a reduced model for fluid flow in fractured
porous media. The liquid phase flow is governed by the Darcy law and, dissolved in the liquid itself, chemical species (solutes) can react and precipitate forming a salt (or an immobile phase that fills the void spaces). Moreover, the
latter can also dissolve to form solutes. The dissolution or
precipitation processes can alter the porosity of the porous media, changing thus
the Darcy velocity of the liquid. As mentioned, we have assumed that in the
porous medium a fracture is present which may dramatically alter the flow properties of the
system and thus requires an adequate model to obtain reliable and accurate
results. What we have proposed is a reduced model that leads to a
hybrid-dimensional framework, where the fracture is one dimensional smaller than
the porous medium itself. New equations have been derived to model the physical
processes in the fracture as well as the coupling conditions between the
fracture itself and the surrounding porous media. The complete set of equations forms a reactive transport model in a fractured porous medium. An extension, which will be part of a future work,  is the introduction of a discrete setting for the efficient solution
of the proposed mathematical model.

%\begin{acknowledgement}
%    If you want to include acknowledgments of assistance and the like at the end of an individual chapter please use the \verb|acknowledgement| environment -- it will automatically render Springer's preferred layout.
%\end{acknowledgement}

% the bibtex style to be used is spmpsci.bst
%\bibliographystyle{spmpsci}
%\bibliography{biblio}

\end{document}

%% file: fracture.pdf_tex
%% Accompanies image file 'fracture.pdf' (pdf, eps, ps)
\begingroup%
  \makeatletter%
  \providecommand\color[2][]{%
    \errmessage{(Inkscape) Color is used for the text in Inkscape, but the package 'color.sty' is not loaded}%
    \renewcommand\color[2][]{}%
  }%
  \providecommand\transparent[1]{%
    \errmessage{(Inkscape) Transparency is used (non-zero) for the text in Inkscape, but the package 'transparent.sty' is not loaded}%
    \renewcommand\transparent[1]{}%
  }%
  \providecommand\rotatebox[2]{#2}%
  \newcommand*\fsize{\dimexpr\f@size pt\relax}%
  \newcommand*\lineheight[1]{\fontsize{\fsize}{#1\fsize}\selectfont}%
  \ifx\svgwidth\undefined%
    \setlength{\unitlength}{154.02376791bp}%
    \ifx\svgscale\undefined%
      \relax%
    \else%
      \setlength{\unitlength}{\unitlength * \real{\svgscale}}%
    \fi%
  \else%
    \setlength{\unitlength}{\svgwidth}%
  \fi%
  \global\let\svgwidth\undefined%
  \global\let\svgscale\undefined%
  \makeatother%
  \begin{picture}(1,0.56888735)%
    \lineheight{1}%
    \setlength\tabcolsep{0pt}%
    \put(0,0){\includegraphics[width=\unitlength,page=1]{fracture.pdf}}%
    \put(0.50450213,0.25390337){\color[rgb]{0,0,0}\makebox(0,0)[lt]{\lineheight{1.25}\smash{\begin{tabular}[t]{l}$\epsilon(t)$\end{tabular}}}}%
    \put(0,0){\includegraphics[width=\unitlength,page=2]{fracture.pdf}}%
    \put(0.09480467,0.29390781){\color[rgb]{0,0,0}\rotatebox{7.7021716}{\makebox(0,0)[lt]{\lineheight{1.25}\smash{\begin{tabular}[t]{l}$\gamma$\end{tabular}}}}}%
    \put(0,0){\includegraphics[width=\unitlength,page=3]{fracture.pdf}}%
    \put(0.69517127,0.44748661){\color[rgb]{0,0,0}\makebox(0,0)[lt]{\lineheight{1.25}\smash{\begin{tabular}[t]{l}$\bm{n}$\end{tabular}}}}%
    \put(0.44083698,0.04424196){\color[rgb]{0,0,0}\makebox(0,0)[lt]{\lineheight{1.25}\smash{\begin{tabular}[t]{l}$\Omega^-$\end{tabular}}}}%
    \put(0.74826292,0.36140652){\color[rgb]{0,0,0}\makebox(0,0)[lt]{\lineheight{1.25}\smash{\begin{tabular}[t]{l}$\Omega_\gamma(t)$\end{tabular}}}}%
    \put(0,0){\includegraphics[width=\unitlength,page=4]{fracture.pdf}}%
    \put(0.18332209,0.48410158){\color[rgb]{0,0,0}\makebox(0,0)[lt]{\lineheight{1.25}\smash{\begin{tabular}[t]{l}$\bm{n}_{\Omega_\gamma}^+(t)$\end{tabular}}}}%
    \put(0,0){\includegraphics[width=\unitlength,page=5]{fracture.pdf}}%
    \put(0.20722959,0.12762587){\color[rgb]{0,0,0}\makebox(0,0)[lt]{\lineheight{1.25}\smash{\begin{tabular}[t]{l}$\bm{n}_{\Omega_\gamma}^-(t)$\end{tabular}}}}%
    \put(0.4561693,0.49800283){\color[rgb]{0,0,0}\makebox(0,0)[lt]{\lineheight{1.25}\smash{\begin{tabular}[t]{l}$\Omega^+$\end{tabular}}}}%
  \end{picture}%
\endgroup%

%% file: fracture_reduced.pdf_tex
%% Accompanies image file 'fracture_reduced.pdf' (pdf, eps, ps)
\begingroup%
  \makeatletter%
  \providecommand\color[2][]{%
    \errmessage{(Inkscape) Color is used for the text in Inkscape, but the package 'color.sty' is not loaded}%
    \renewcommand\color[2][]{}%
  }%
  \providecommand\transparent[1]{%
    \errmessage{(Inkscape) Transparency is used (non-zero) for the text in Inkscape, but the package 'transparent.sty' is not loaded}%
    \renewcommand\transparent[1]{}%
  }%
  \providecommand\rotatebox[2]{#2}%
  \newcommand*\fsize{\dimexpr\f@size pt\relax}%
  \newcommand*\lineheight[1]{\fontsize{\fsize}{#1\fsize}\selectfont}%
  \ifx\svgwidth\undefined%
    \setlength{\unitlength}{154.02376791bp}%
    \ifx\svgscale\undefined%
      \relax%
    \else%
      \setlength{\unitlength}{\unitlength * \real{\svgscale}}%
    \fi%
  \else%
    \setlength{\unitlength}{\svgwidth}%
  \fi%
  \global\let\svgwidth\undefined%
  \global\let\svgscale\undefined%
  \makeatother%
  \begin{picture}(1,0.56888735)%
    \lineheight{1}%
    \setlength\tabcolsep{0pt}%
    \put(0,0){\includegraphics[width=\unitlength,page=1]{fracture_reduced.pdf}}%
    \put(0.09480467,0.29390781){\color[rgb]{0,0,0}\rotatebox{7.7021716}{\makebox(0,0)[lt]{\lineheight{1.25}\smash{\begin{tabular}[t]{l}$\gamma$\end{tabular}}}}}%
    \put(0,0){\includegraphics[width=\unitlength,page=2]{fracture_reduced.pdf}}%
    \put(0.69517127,0.44748661){\color[rgb]{0,0,0}\makebox(0,0)[lt]{\lineheight{1.25}\smash{\begin{tabular}[t]{l}$\bm{n}$\end{tabular}}}}%
    \put(0.44083698,0.04424196){\color[rgb]{0,0,0}\makebox(0,0)[lt]{\lineheight{1.25}\smash{\begin{tabular}[t]{l}$\Omega^-$\end{tabular}}}}%
    \put(0,0){\includegraphics[width=\unitlength,page=3]{fracture_reduced.pdf}}%
    \put(0.4561693,0.49800283){\color[rgb]{0,0,0}\makebox(0,0)[lt]{\lineheight{1.25}\smash{\begin{tabular}[t]{l}$\Omega^+$\end{tabular}}}}%
  \end{picture}%
\endgroup%

%% file: fracture_omega.pdf_tex
%% Accompanies image file 'fracture_omega.pdf' (pdf, eps, ps)
\begingroup%
  \makeatletter%
  \providecommand\color[2][]{%
    \errmessage{(Inkscape) Color is used for the text in Inkscape, but the package 'color.sty' is not loaded}%
    \renewcommand\color[2][]{}%
  }%
  \providecommand\transparent[1]{%
    \errmessage{(Inkscape) Transparency is used (non-zero) for the text in Inkscape, but the package 'transparent.sty' is not loaded}%
    \renewcommand\transparent[1]{}%
  }%
  \providecommand\rotatebox[2]{#2}%
  \newcommand*\fsize{\dimexpr\f@size pt\relax}%
  \newcommand*\lineheight[1]{\fontsize{\fsize}{#1\fsize}\selectfont}%
  \ifx\svgwidth\undefined%
    \setlength{\unitlength}{154.02376791bp}%
    \ifx\svgscale\undefined%
      \relax%
    \else%
      \setlength{\unitlength}{\unitlength * \real{\svgscale}}%
    \fi%
  \else%
    \setlength{\unitlength}{\svgwidth}%
  \fi%
  \global\let\svgwidth\undefined%
  \global\let\svgscale\undefined%
  \makeatother%
  \begin{picture}(1,0.56888735)%
    \lineheight{1}%
    \setlength\tabcolsep{0pt}%
    \put(0,0){\includegraphics[width=\unitlength,page=1]{fracture_omega.pdf}}%
    \put(0.76797044,0.25651193){\color[rgb]{0,0,0}\makebox(0,0)[lt]{\lineheight{1.25}\smash{\begin{tabular}[t]{l}$\epsilon$\end{tabular}}}}%
    \put(0,0){\includegraphics[width=\unitlength,page=2]{fracture_omega.pdf}}%
    \put(0.02089468,0.28086477){\color[rgb]{0,0,0}\rotatebox{7.7021716}{\makebox(0,0)[lt]{\lineheight{1.25}\smash{\begin{tabular}[t]{l}$\gamma$\end{tabular}}}}}%
    \put(0,0){\includegraphics[width=\unitlength,page=3]{fracture_omega.pdf}}%
    \put(0.69517127,0.44748661){\color[rgb]{0,0,0}\makebox(0,0)[lt]{\lineheight{1.25}\smash{\begin{tabular}[t]{l}$\bm{n}$\end{tabular}}}}%
    \put(0.44083698,0.04424196){\color[rgb]{0,0,0}\makebox(0,0)[lt]{\lineheight{1.25}\smash{\begin{tabular}[t]{l}$\Omega^-$\end{tabular}}}}%
    \put(0.85869335,0.35618941){\color[rgb]{0,0,0}\makebox(0,0)[lt]{\lineheight{1.25}\smash{\begin{tabular}[t]{l}$\Omega_\gamma$\end{tabular}}}}%
    \put(0,0){\includegraphics[width=\unitlength,page=4]{fracture_omega.pdf}}%
    \put(0.42939966,0.1588968){\color[rgb]{0,0,0}\makebox(0,0)[lt]{\lineheight{1.25}\smash{\begin{tabular}[t]{l}$\bm{n}_{\omega}$\end{tabular}}}}%
    \put(0,0){\includegraphics[width=\unitlength,page=5]{fracture_omega.pdf}}%
    \put(0.4561693,0.49800283){\color[rgb]{0,0,0}\makebox(0,0)[lt]{\lineheight{1.25}\smash{\begin{tabular}[t]{l}$\Omega^+$\end{tabular}}}}%
    \put(0,0){\includegraphics[width=\unitlength,page=6]{fracture_omega.pdf}}%
    \put(0.36254381,0.25941703){\color[rgb]{0,0,0}\makebox(0,0)[lt]{\lineheight{1.25}\smash{\begin{tabular}[t]{l}$\omega$\end{tabular}}}}%
    \put(0.37915106,0.45518618){\color[rgb]{0,0,0}\makebox(0,0)[lt]{\lineheight{1.25}\smash{\begin{tabular}[t]{l}$\bm{n}_{\omega}$\end{tabular}}}}%
    \put(0,0){\includegraphics[width=\unitlength,page=7]{fracture_omega.pdf}}%
    \put(0.16020258,0.21536898){\color[rgb]{0,0,0}\makebox(0,0)[lt]{\lineheight{1.25}\smash{\begin{tabular}[t]{l}$\bm{n}_{\omega}$\end{tabular}}}}%
    \put(0.22364364,0.3161105){\color[rgb]{0,0,0}\makebox(0,0)[lt]{\lineheight{1.25}\smash{\begin{tabular}[t]{l}$l_0$\end{tabular}}}}%
    \put(0.6233269,0.35101611){\color[rgb]{0,0,0}\makebox(0,0)[lt]{\lineheight{1.25}\smash{\begin{tabular}[t]{l}$l_1$\end{tabular}}}}%
  \end{picture}%
\endgroup%

%% file: fracture_aperture.pdf_tex
%% Accompanies image file 'fracture_aperture.pdf' (pdf, eps, ps)
\begingroup%
  \makeatletter%
  \providecommand\color[2][]{%
    \errmessage{(Inkscape) Color is used for the text in Inkscape, but the package 'color.sty' is not loaded}%
    \renewcommand\color[2][]{}%
  }%
  \providecommand\transparent[1]{%
    \errmessage{(Inkscape) Transparency is used (non-zero) for the text in Inkscape, but the package 'transparent.sty' is not loaded}%
    \renewcommand\transparent[1]{}%
  }%
  \providecommand\rotatebox[2]{#2}%
  \newcommand*\fsize{\dimexpr\f@size pt\relax}%
  \newcommand*\lineheight[1]{\fontsize{\fsize}{#1\fsize}\selectfont}%
  \ifx\svgwidth\undefined%
    \setlength{\unitlength}{154.02376791bp}%
    \ifx\svgscale\undefined%
      \relax%
    \else%
      \setlength{\unitlength}{\unitlength * \real{\svgscale}}%
    \fi%
  \else%
    \setlength{\unitlength}{\svgwidth}%
  \fi%
  \global\let\svgwidth\undefined%
  \global\let\svgscale\undefined%
  \makeatother%
  \begin{picture}(1,0.56888735)%
    \lineheight{1}%
    \setlength\tabcolsep{0pt}%
    \put(0,0){\includegraphics[width=\unitlength,page=1]{fracture_aperture.pdf}}%
    \put(0.50450213,0.25390337){\color[rgb]{0,0,0}\makebox(0,0)[lt]{\lineheight{1.25}\smash{\begin{tabular}[t]{l}$\epsilon(t)$\end{tabular}}}}%
    \put(0,0){\includegraphics[width=\unitlength,page=2]{fracture_aperture.pdf}}%
    \put(0.44083698,0.04424196){\color[rgb]{0,0,0}\makebox(0,0)[lt]{\lineheight{1.25}\smash{\begin{tabular}[t]{l}$\Omega^-$\end{tabular}}}}%
    \put(0.74826292,0.36140652){\color[rgb]{0,0,0}\makebox(0,0)[lt]{\lineheight{1.25}\smash{\begin{tabular}[t]{l}$\Omega_\gamma(t)$\end{tabular}}}}%
    \put(0,0){\includegraphics[width=\unitlength,page=3]{fracture_aperture.pdf}}%
    \put(0.4561693,0.49800283){\color[rgb]{0,0,0}\makebox(0,0)[lt]{\lineheight{1.25}\smash{\begin{tabular}[t]{l}$\Omega^+$\end{tabular}}}}%
    \put(0,0){\includegraphics[width=\unitlength,page=4]{fracture_aperture.pdf}}%
  \end{picture}%
\endgroup%

%% file: alessio_fumagalli.bbl
\begin{thebibliography}{10}
\providecommand{\url}[1]{{#1}}
\providecommand{\urlprefix}{URL }
\expandafter\ifx\csname urlstyle\endcsname\relax
  \providecommand{\doi}[1]{DOI~\discretionary{}{}{}#1}\else
  \providecommand{\doi}{DOI~\discretionary{}{}{}\begingroup
  \urlstyle{rm}\Url}\fi

\bibitem{Agosti2015}
Agosti, A., Formaggia, L., Scotti, A.: Analysis of a model for precipitation
  and dissolution coupled with a darcy flux.
\newblock Journal of Mathematical Analysis and Applications \textbf{431}(2),
  752--781 (2015).
\newblock \doi{https://doi.org/10.1016/j.jmaa.2015.06.003}.
\newblock
  \urlprefix\url{http://www.sciencedirect.com/science/article/pii/S0022247X15005466}

\bibitem{Alboin2002}
Alboin, C., Jaffr{\'e}, J., Roberts, J.E., Serres, C.: Modeling fractures as
  interfaces for flow and transport in porous media.
\newblock In: Fluid flow and transport in porous media: mathematical and
  numerical treatment ({S}outh {H}adley, {MA}, 2001), \emph{Contemp. Math.},
  vol. 295, pp. 13--24. Amer. Math. Soc., Providence, RI (2002)

\bibitem{Ambartsumyan2019}
Ambartsumyan, I., Khattatov, E., Nguyen, T., Yotov, I.: Flow and transport in
  fractured poroelastic media.
\newblock GEM - International Journal on Geomathematics \textbf{10}(1), 11
  (2019).
\newblock \doi{10.1007/s13137-019-0119-5}.
\newblock \urlprefix\url{https://doi.org/10.1007/s13137-019-0119-5}

\bibitem{Antonietti}
Antonietti, P.F., Formaggia, L., Scotti, A., Verani, M., Verzotti, N.: Mimetic
  finite difference approximation of flows in fractured porous media.
\newblock ESAIM: M2AN \textbf{50}(3), 809--832 (2016).
\newblock \doi{10.1051/m2an/2015087}.
\newblock \urlprefix\url{http://dx.doi.org/10.1051/m2an/2015087}

\bibitem{Berge2019}
Berge, R., Berre, I., Keilegavlen, E., Wohlmuth, B.: Finite volume
  discretization for poroelastic media with fractures modeled by contact
  mechanics.
\newblock Tech. rep., arXiv:1904.11916 [math.NA] (2019).
\newblock \urlprefix\url{https://arxiv.org/abs/1904.11916}

\bibitem{Berre2019b}
Berre, I., Doster, F., Keilegavlen, E.: Flow in fractured porous media: A
  review of conceptual models and discretization approaches.
\newblock Transport in Porous Media \textbf{130}(1), 215--236 (2019).
\newblock \doi{10.1007/s11242-018-1171-6}.
\newblock \urlprefix\url{https://doi.org/10.1007/s11242-018-1171-6}

\bibitem{Boon2018}
Boon, W.M., Nordbotten, J.M., Yotov, I.: Robust discretization of flow in
  fractured porous media.
\newblock SIAM Journal on Numerical Analysis \textbf{56}(4), 2203--2233 (2018).
\newblock \doi{10.1137/17M1139102}.
\newblock \urlprefix\url{https://doi.org/10.1137/17M1139102}

\bibitem{Bukac2017}
Bukac, M., Yotov, I., Zunino, P.: Dimensional model reduction for flow through
  fractures in poroelastic media.
\newblock ESAIM: Mathematical Modelling and Numerical Analysis \textbf{51}(4),
  1429--1471 (2017).
\newblock \doi{10.1051/m2an/2016069}.
\newblock \urlprefix\url{https://doi.org/10.1051/m2an/2016069}

\bibitem{Chave2018}
Chave, F., Di~Pietro, D.A., Formaggia, L.: A hybrid high-order method for darcy
  flows in fractured porous media.
\newblock SIAM Journal on Scientific Computing \textbf{40}(2), A1063--A1094
  (2018).
\newblock \doi{10.1137/17M1119500}.
\newblock \urlprefix\url{https://doi.org/10.1137/17M1119500}

\bibitem{Chave2019}
Chave, F., Di~Pietro, D.A., Formaggia, L.: A hybrid high-order method for
  passive transport in fractured porous media.
\newblock GEM - International Journal on Geomathematics \textbf{10}(1), 12
  (2019).
\newblock \doi{10.1007/s13137-019-0114-x}.
\newblock \urlprefix\url{https://doi.org/10.1007/s13137-019-0114-x}

\bibitem{DAngelo2011}
D'Angelo, C., Scotti, A.: A mixed finite element method for {D}arcy flow in
  fractured porous media with non-matching grids.
\newblock Mathematical {M}odelling and {N}umerical {A}nalysis \textbf{46}(02),
  465--489 (2012).
\newblock \doi{10.1051/m2an/2011148}

\bibitem{Duijn2004}
van Duijn, C., Pop, I.S.: Crystal dissolution and precipitation in porous media
  : pore scale analysis.
\newblock Journal f{\"u}r die reine und angewandte Mathematik (Crelle's
  Journal) \textbf{577}, 171--211 (2004).
\newblock \doi{10.1515/crll.2004.2004.577.171}

\bibitem{Elyes2015}
Elyes, A., J{\'e}r{\^o}me, J., Roberts, J.E.: A 3-{D} reduced fracture model
  for two-phase flow in porous media with a global pressure formulation.
\newblock In: {MAMERN VI}. Pau, France (2015).
\newblock \urlprefix\url{https://hal.inria.fr/hal-01119986}

\bibitem{Emmanuel2005}
Emmanuel, S., Berkowitz, B.: Mixing-induced precipitation and porosity
  evolution in porous media.
\newblock Advances in Water Resources \textbf{28}(4), 337--344 (2005).
\newblock \doi{https://doi.org/10.1016/j.advwatres.2004.11.010}.
\newblock
  \urlprefix\url{http://www.sciencedirect.com/science/article/pii/S030917080400212X}

\bibitem{Faille2014a}
Faille, I., Fumagalli, A., Jaffr{\'e}, J., Roberts, J.E.: Model reduction and
  discretization using hybrid finite volumes of flow in porous media containing
  faults.
\newblock Computational Geosciences \textbf{20}(2), 317--339 (2016).
\newblock \doi{10.1007/s10596-016-9558-3}.
\newblock
  \urlprefix\url{https://link.springer.com/article/10.1007/s10596-016-9558-3}

\bibitem{Formaggia2012}
Formaggia, L., Fumagalli, A., Scotti, A., Ruffo, P.: A reduced model for
  {D}arcy's problem in networks of fractures.
\newblock {ESAIM}: {M}athematical {M}odelling and {N}umerical {A}nalysis
  \textbf{48}, 1089--1116 (2014).
\newblock \doi{10.1051/m2an/2013132}.
\newblock
  \urlprefix\url{https://www.esaim-m2an.org/articles/m2an/abs/2014/04/m2an130132/m2an130132.html}

\bibitem{Fumagalli2017a}
Fumagalli, A., Keilegavlen, E.: Dual virtual element methods for discrete
  fracture matrix models.
\newblock Oil \& Gas Science and Technology - Revue d'IFP Energies nouvelles
  \textbf{74}(41), 1--17 (2019).
\newblock \doi{doi.org/10.2516/ogst/2019008}

\bibitem{Fumagalli2012d}
Fumagalli, A., Scotti, A.: A numerical method for two-phase flow in fractured
  porous media with non-matching grids.
\newblock Advances in Water Resources \textbf{62, Part C}(0), 454--464 (2013).
\newblock \doi{10.1016/j.advwatres.2013.04.001}.
\newblock
  \urlprefix\url{https://www.sciencedirect.com/science/article/pii/S0309170813000523}.
\newblock Computational Methods in Geologic CO2 Sequestration

\bibitem{Fumagalli2012a}
Fumagalli, A., Scotti, A.: A {R}educed {M}odel for {F}low and {T}ransport in
  {F}ractured {P}orous {M}edia with {N}on-matching {G}rids.
\newblock In: A.~Cangiani, R.L. Davidchack, E.~Georgoulis, A.N. Gorban,
  J.~Levesley, M.V. Tretyakov (eds.) Numerical Mathematics and Advanced
  Applications 2011, pp. 499--507. Springer Berlin Heidelberg (2013).
\newblock \doi{10.1007/978-3-642-33134-3\_53}

\bibitem{Ganis2014}
Ganis, B., Girault, V., Mear, M., Singh, G., Wheeler, M.: Modeling fractures in
  a poro-elastic medium.
\newblock Oil \& Gas Science and Technology--Revue d'IFP Energies nouvelles
  \textbf{69}(4), 515--528 (2014)

\bibitem{Hommel2018}
Hommel, J., Coltman, E., Class, H.: Porosity-permeability relations for
  evolving pore space: A review with a focus on (bio-)geochemically altered
  porous media.
\newblock Transport in Porous Media \textbf{124}(2), 589--629 (2018).
\newblock \doi{10.1007/s11242-018-1086-2}.
\newblock \urlprefix\url{https://doi.org/10.1007/s11242-018-1086-2}

\bibitem{Hornung1994}
Hornung, U., J\"ager, W., Mikeli\'c, A.: Reactive transport through an array of
  cells with semi-permeable membranes.
\newblock ESAIM: Mathematical Modelling and Numerical Analysis - Mod\'elisation
  Math\'ematique et Analyse Num\'erique \textbf{28}(1), 59--94 (1994)

\bibitem{Jaffre2011}
Jaffr{\'e}, J., Mnejja, M., Roberts, J.E.: A discrete fracture model for
  two-phase flow with matrix-fracture interaction.
\newblock Procedia Computer Science \textbf{4}, 967--973 (2011).
\newblock \doi{10.1016/j.procs.2011.04.102}.
\newblock
  \urlprefix\url{http://www.sciencedirect.com/science/article/pii/S1877050911001608}

\bibitem{Katz2011}
Katz, G.E., Berkowitz, B., Guadagnini, A., Saaltink, M.W.: Experimental and
  modeling investigation of multicomponent reactive transport in porous media.
\newblock Journal of Contaminant Hydrology \textbf{120--121}, 27--44 (2011).
\newblock \doi{https://doi.org/10.1016/j.jconhyd.2009.11.002}.
\newblock
  \urlprefix\url{http://www.sciencedirect.com/science/article/pii/S0169772209001582}.
\newblock Reactive Transport in the Subsurface: Mixing, Spreading and Reaction
  in Heterogeneous Media

\bibitem{Knabner1995}
Knabner, P., van Duijn, C., Hengst, S.: An analysis of crystal dissolution
  fronts in flows through porous media. part 1: Compatible boundary conditions.
\newblock Advances in Water Resources \textbf{18}(3), 171--185 (1995).
\newblock \doi{https://doi.org/10.1016/0309-1708(95)00005-4}.
\newblock
  \urlprefix\url{http://www.sciencedirect.com/science/article/pii/0309170895000054}

\bibitem{Kumar2011}
Kumar, K., van Noorden, T.L., Pop, I.S.: Effective dispersion equations for
  reactive flows involving free boundaries at the microscale.
\newblock Multiscale Modeling \& Simulation \textbf{9}(1), 29--58 (2011).
\newblock \doi{10.1137/100804553}.
\newblock \urlprefix\url{https://doi.org/10.1137/100804553}

\bibitem{Martin2005}
Martin, V., Jaffr{\'e}, J., Roberts, J.E.: Modeling {F}ractures and {B}arriers
  as {I}nterfaces for {F}low in {P}orous {M}edia.
\newblock SIAM J. Sci. Comput. \textbf{26}(5), 1667--1691 (2005).
\newblock \doi{10.1137/S1064827503429363}

\bibitem{Noorden2009a}
van Noorden, T.: Crystal precipitation and dissolution in a thin strip.
\newblock European Journal of Applied Mathematics \textbf{20}(1), 69--91
  (2009).
\newblock \doi{10.1017/S0956792508007651}

\bibitem{Noorden2009}
van Noorden, T.L.: Crystal precipitation and dissolution in a porous medium:
  Effective equations and numerical experiments.
\newblock Multiscale Modeling \& Simulation \textbf{7}(3), 1220--1236 (2009).
\newblock \doi{10.1137/080722096}.
\newblock \urlprefix\url{https://doi.org/10.1137/080722096}

\bibitem{Ray2019}
Ray, N., Oberlander, J., Frolkovic, P.: Numerical investigation of a fully
  coupled micro-macro model for mineral dissolution and precipitation.
\newblock Computational Geosciences \textbf{23}(5), 1173--1192 (2019).
\newblock \doi{10.1007/s10596-019-09876-x}.
\newblock \urlprefix\url{https://doi.org/10.1007/s10596-019-09876-x}

\bibitem{Sandve2012}
Sandve, T.H., Berre, I., Nordbotten, J.M.: An efficient multi-point flux
  approximation method for {D}iscrete {F}racture-{M}atrix simulations.
\newblock Journal of Computational Physics \textbf{231}(9), 3784--3800 (2012).
\newblock \doi{https://doi.org/10.1016/j.jcp.2012.01.023}.
\newblock
  \urlprefix\url{http://www.sciencedirect.com/science/article/pii/S0021999112000447}

\bibitem{Schwenck2015}
Schwenck, N., Flemisch, B., Helmig, R., Wohlmuth, B.: Dimensionally reduced
  flow models in fractured porous media: crossings and boundaries.
\newblock Computational Geosciences \textbf{19}(6), 1219--1230 (2015).
\newblock \doi{10.1007/s10596-015-9536-1}.
\newblock \urlprefix\url{http://dx.doi.org/10.1007/s10596-015-9536-1}

\bibitem{Scotti2017}
Scotti, A., Formaggia, L., Sottocasa, F.: Analysis of a mimetic finite
  difference approximation of flows in fractured porous media.
\newblock ESAIM: M2AN  (2017).
\newblock \doi{doi.org/10.1051/m2an/2017028}

\bibitem{Stefansson2018}
Stefansson, I., Berre, I., Keilegavlen, E.: Finite-volume discretisations for
  flow in fractured porous media.
\newblock Transport in Porous Media \textbf{124}(2), 439--462 (2018).
\newblock \doi{10.1007/s11242-018-1077-3}.
\newblock \urlprefix\url{https://doi.org/10.1007/s11242-018-1077-3}

\bibitem{Ucar2018}
Ucar, E., Keilegavlen, E., Berre, I., Nordbotten, J.M.: A finite-volume
  discretization for deformation of fractured media.
\newblock Computational Geosciences \textbf{22}(4), 993--1007 (2018).
\newblock \doi{10.1007/s10596-018-9734-8}.
\newblock \urlprefix\url{https://doi.org/10.1007/s10596-018-9734-8}

\end{thebibliography}
